\documentstyle[makeidx,draft,12pt]{article}

\newtheorem{th}{Theorem}[section]

\newtheorem{prop}[th]{Proposition}

\newtheorem{defn}[th]{Definition}
\newenvironment{defn-new}{\begin{defn} \em}{\end{defn}}
\newtheorem{rem}[th]{Remark}
\newenvironment{rem-new}{\begin{rem} \em}{\end{rem}}
\newtheorem{ex}[th]{Example}
\newenvironment{ex-new}{\begin{ex} \em}{\end{ex}}

\newenvironment{notation-new}{\begin{rem} \em}{\end{rem}}

\newenvironment{agr-new}{\begin{rem} \em}{\end{rem}}

\makeatletter \@addtoreset{equation}{section} \makeatother

\makeatletter \@addtoreset{figure}{section} \makeatother

\begin{document}

\begin{center}
{\Large On Cosymplectic Conformal Connections}

{\large Punam Gupta}

({\small Dedicated to K.Yano)}
\end{center}

\noindent {\bf Abstract.} The aim of this paper is to introduce a
cosymplectic analouge of conformal connection in a cosymplectic manifold and
proved that if cosymplectic manifold $M$ admits a cosymplectic conformal
connection which is of zero curvature, then the Bochner curvature tensor of $%
M$ vanishes.\medskip

\noindent {\bf Mathematics Subject Classification.}53C05, 53C18, 53C20
53D05. \medskip

\noindent {\bf Keywords. }K\"{a}hler manifold, Sasakian manifold,
Cosymplectic manifold.

\section{Introduction}

In 1958, Libermann \cite{Liber, Liber1} introduced the cosymplectic
manifolds. He defined it as :The pair $(\eta ,\omega )$ of a $1$-form $\eta $
and a $2$-form $\omega $ such that $\eta \wedge \omega ^{n}$ is a volume
form on odd-dimensional manifold $M^{2m+1}$, defines an almost cosymplectic
structure on $M$. If both $\eta $ and $\omega $ are closed, the structure is
said to be cosymplectic. So, a manifold endowed with a cosymplectic
structure $(M,\eta ,\omega )$ is called a cosymplectic manifold.

After few years, in 1967, Blair \cite{Blair1} used the term
\textquotedblleft cosymplectic\textquotedblright\ to define manifold with
almost contact metric structure satisfying a normality condition, although
related to the Libermann's definition.

In differential geometry, one of the most important fields is to study
smooth maps which preserve certain geometric properties. In those
transfprmations, the conformal transformation is very interesting, where
only angles are preserved both in magnitude and orientation but not
necessarily distances. Stereographic projection is the simplest example of
conformal transformation. It is belived that, in 1569, the property of
conformal transformations is first used by Gerardus Mercator to produce the
famous Mercator's world map (the first angle-preserving (or conformal) world
map). For more details on conformal transformations, see \cite{Bob}.

\begin{defn-new}
A diffeomorphism $f:(M,g)\rightarrow (N,\tilde{g})$ is said to be conformal
mapping \cite{Eisenhart} if
\begin{equation}
\tilde{g}=e^{2p}g,  \label{eq-c}
\end{equation}%
where $p$ is a function on $M$. If $p$ is constant, then conformal mapping
is called homothetic mapping.
\end{defn-new}

Let $\nabla $ be a linear connection in an $n$-dimensioanl manifold $M$. The
torsion tensor $T_{ji}^{h}$ of $\nabla $ is
\[
T_{ji}^{h}=\Gamma _{ji}^{h}-\Gamma _{ij}^{h}.
\]%
The connection $\nabla $ is symmetric if the torsion tensor $T_{ji}^{h}$
vanishes, otherwise it is non-symmetric. The connection $\nabla $ is metric
if there is a Riemannian metric $g_{ji}$ in $M$ such that $\nabla
_{k}g_{ji}=0$, otherwise it is non-metric. It is well known that a linear
connection is symmetric and metric if and only if it is the Levi-Civita
connection.

Let $D$ be covariant differentiation operator with respect to Christofell
symbols $\Gamma _{ji}^{h}$ formed with $\tilde{g}$ corresponding to the
conformal change (\ref{eq-c}), then we have
\begin{equation}
\Gamma _{ji}^{h}=\left\{
\begin{array}{c}
h \\
ji%
\end{array}%
\right\} +\delta _{j}^{h}p_{i}+\delta _{i}^{h}p_{j}-g_{ji}p^{h},
\label{eq-con}
\end{equation}%
where $p_{i}$ is the gradient of $p$, that is, $p_{i}=\partial _{i}p$ and $%
p^{h}=p_{t}g^{th}$, $g^{th}$ is the contravariant components of the metric
tensor $g_{th}$. This affine connection is known as conformal connection.

So, we have
\[
D_{k}(e^{2p}g_{ji})=0.
\]%
Let $R$ be the curvature tensor of the Riemannian metric $\tilde{g}$. Then
\[
R_{kji}^{h}=K_{kji}^{h}+\delta _{k}^{h}p_{ji}-\delta
_{j}^{h}p_{ki}+p_{k}^{h}g_{ji}-p_{j}^{h}g_{ki},
\]%
where $K$ be the curvature tensor of the Riemannian metric $g$ and
\begin{equation}
p_{ji}=\nabla _{j}p_{i}-p_{j}p_{i}+\frac{1}{2}\lambda p_{t}p^{t}g_{ji}.
\end{equation}%
If the Riemannian metric $g$ is conformally equivalent to the Riemannian
metric $\tilde{g}$, which is locally Euclidean, then the Riemannian manifold
with the metric $g$ is said to be conformally flat \cite{Yano}. For which,
we are stating the well known theorem given by Weyl \cite{Weyl}.

\begin{th}
A necessary and sufficient condition for a Riemannian manifold to be
conformally flat is that Weyl conformal curvature tensor $C=0$ for $n>3$, or
$C_{kji}^{h}=0$, where
\[
C_{kji}^{h}=K_{kji}^{h}+\delta _{k}^{h}C_{ji}-\delta
_{j}^{h}C_{ki}+C_{k}^{h}g_{ji}-C_{j}^{h}g_{ki},
\]%
and
\[
C_{ji}=p_{ji}.
\]
\end{th}

In 1975, K. Yano \cite{Yano-75} studied a complex analogue of conformal
connection and proved the following remarkable theorems:

\begin{th}
{\rm \cite{Yano-75}} In a K\"{a}hlerian manifold with Hermitian metric
tensor $g_{ji}$ and almost complex structure tensor $F_{i}^{h}$, the affine
connection $D$ with components $\Gamma _{ji}^{h}$ which satisfies%
\[
D_{k}(e^{2p}g_{ji})=0,
\]%
\[
D_{k}(e^{2p}\varphi _{ji})=0,
\]%
\[
\Gamma _{ji}^{h}-\Gamma _{ij}^{h}=-2\varphi _{ji}q^{h},
\]%
where $p$ is a scalar function, $q^{h}$ a vector field and $\varphi
_{ji}=\varphi _{j}^{t}g_{ti}$, is given by%
\[
\Gamma _{ji}^{h}=\left\{
\begin{array}{c}
h \\
ji%
\end{array}%
\right\} +\delta _{j}^{h}p_{i}+\delta _{i}^{h}p_{j}-g_{ji}p^{h}+\varphi
_{j}^{h}q_{i}+\varphi _{i}^{h}q_{j}-\varphi _{ji}q^{h},
\]%
where $p_{i}$ is the gradient of $p$, that is, $p_{i}=\partial _{i}p$ and $%
p^{h}=p_{t}g^{th},\quad q_{i}=-p_{t}\varphi _{i}^{t},\quad q^{h}=q_{t}g^{th}$%
.
\end{th}

He called such an affine connection a complex conformal connection in a K%
\"{a}hlerian manifold.

\begin{th}
If in a real $n$-dimensional Kaehlerian manifold, $(n\geq 4)$, there exists
a scalar function $p$ such that the complex conformal connection is of zero
curvature, then the Bochner curvature tensor of the manifold
\begin{eqnarray}
B_{kji}^{h} &=&K_{kji}^{h}+\delta _{k}^{h}L_{ji}-\delta
_{j}^{h}L_{ki}+L_{k}^{h}g_{ji}-L_{j}^{h}g_{ki}  \nonumber \\
&&+F_{k}^{h}M_{ji}-F_{j}^{h}M_{ki}+M_{k}^{h}F_{ji}-M_{j}^{h}F_{ki}  \nonumber
\\
&&-2(M_{kj}F_{i}^{h}+F_{kj}M_{i}^{h})  \nonumber
\end{eqnarray}%
vanishes, where
\[
L_{ji}=-\frac{1}{n+4}K_{ji}+\frac{1}{2(n+2)(n+4)}Kg_{ji},\quad
L_{k}^{h}=L_{ki}g^{th},
\]%
\[
M_{ji}=-L_{jt}F_{i}^{t},\quad M_{k}^{h}=M_{kt}g^{th}.
\]
\end{th}

After that Yano \cite{Yano-77.} studied a contact analouge of conformal
connection and proved the following results:

\begin{th}
In a Sasakian manifold with structure tensors $(\varphi _{i}^{h},\eta
_{i},g_{ji})$, the affine connection $D$ with components $\Gamma _{ji}^{h}$
which satisfies%
\[
D_{k}(e^{2p}g_{ji})=2e^{2p}p_{k}\eta _{j}\eta _{i},\quad D_{j}\varphi
_{i}^{h}=0,\quad D_{j}\eta ^{h}=0
\]%
\[
\Gamma _{ji}^{h}-\Gamma _{ij}^{h}=-2\varphi _{ji}u^{h},
\]%
where $p$ is a scalar function and $u^{h}$ a vector field, is given by%
\begin{eqnarray*}
\Gamma _{ji}^{h} &=&\left\{
\begin{array}{c}
h \\
ji%
\end{array}%
\right\} +\left( \delta _{j}^{h}-\eta _{j}\eta ^{h}\right) p_{i}+\left(
\delta _{i}^{h}-\eta _{i}\eta ^{h}\right) p_{j}-\left( g_{ji}-\eta _{j}\eta
_{i}\right) p^{h} \\
&&+\varphi _{j}^{h}(q_{i}-\eta _{i})+\varphi _{i}^{h}(q_{j}-\eta
_{j})-\varphi _{ji}(q^{h}-\eta ^{h}),
\end{eqnarray*}%
where%
\[
p_{i}=\partial _{i}p,\quad p^{h}=p_{t}g^{th},\quad q_{i}=-p_{t}\varphi
_{i}^{t},\quad q^{h}=q_{t}g^{th},
\]%
and $p$ satisfies
\[
{\cal L}p=0.
\]
\end{th}

He called such an affine connection a contact conformal connection in a
Sasakian manifold.

\begin{th}
If, in a $n(=2m+1)$-dimensional Sasakian manifold, $(n\geq 4)$, there exists
a scalar function $p$ such that the contact conformal connection is of zero
curvature, then the contact Bochner curvature tensor of the manifold
\begin{eqnarray*}
B_{kji}^{h} &=&K_{kji}^{h}+(\delta _{k}^{h}-\eta _{k}\eta
^{h})L_{ji}-(\delta _{j}^{h}-\eta _{j}\eta ^{h})L_{ki} \\
&&+L_{k}^{h}(g_{ji}-\eta _{j}\eta _{i})-L_{j}^{h}(g_{ki}-\eta _{k}\eta _{i})
\\
&&+\varphi _{k}^{h}M_{ji}-\varphi _{j}^{h}M_{ki}+M_{k}^{h}\varphi
_{ji}-M_{j}^{h}\varphi _{ki} \\
&&-2(M_{kj}\varphi _{i}^{h}+\varphi _{kj}M_{i}^{h})+(\varphi _{k}^{h}\varphi
_{ji}-\varphi _{j}^{h}\varphi _{ki}-2\varphi _{kj}\varphi _{i}^{h}),
\end{eqnarray*}%
vanishes, where
\[
L_{ji}=-\frac{1}{2(m+2)}\left( K_{ji}+(L+3)g_{ji}-(L-1)\eta _{j}\eta
_{i}\right) ,
\]%
\[
L_{k}^{h}=L_{kt}g^{th},\quad M_{ji}=-L_{jt}\varphi _{i}^{t},\quad
M_{k}^{h}=M_{kt}g^{th},\quad L=g^{ji}L_{ji},
\]
\end{th}

The main purpose of the present paper is to introduce a cosymplectic
analouge of conformal connection in a cosymplectic manifold and study its
properties. In Section \ref{sect-Sasakian}, we state some of preliminaries
on cosymplectic manifold and on the cosymplectic Bochner curvature tensor.
In section \ref{sect-Co}, we introduce cosymplectic conformal connection and
in the last section \ref{sect-curv}, we prove that if cosymplectic manifold $%
M$ admits a cosymplectic conformal connection which is of zero curvature,
then the Bochner curvature tensor of $M$ vanishes.

\section{Preliminaries\label{sect-Sasakian}}

\subsection{Cosymplectic manifolds}

Let $M$ be a $(2m+1)$-dimensional differentiable manifold of class $%
C^{\infty }$ covered by a system of coordinate neighborhoods $\left\{
U,y^{h}\right\} $. Let $M$ admit an almost contact structure $(\varphi
_{i}^{h},\xi ^{h},\eta _{i},g_{ij})$ of a tensor field $\varphi _{i}^{h}$ of
type $(1,1)$, a vector field $\xi ^{h}$, a $1$-form $\eta _{i}$ and a
Riemannian metric $g_{ij}$ satisfying
\begin{equation}
\varphi _{j}^{i}\varphi _{i}^{h}=-\delta _{j}^{h}+\eta _{j}\xi ^{h},\quad
\varphi _{i}^{h}\xi ^{i}=0,\quad \eta _{i}\varphi _{j}^{i}=0,\quad \eta
_{i}\xi ^{i}=1,  \label{eq-1}
\end{equation}%
where the indices $h,i,j,k,...$\ $\in \left\{ 1,2,\ldots ,2n+1\right\} $. A
manifold $M$ with an almost contact structure is known as an {\em almost
contact manifold} \cite{Blair,Sasaki-65}.

If the set $(\varphi _{i}^{h},\xi ^{h},\eta _{i})$ satisfies
\begin{equation}
N_{ji}^{h}+(\partial _{j}\eta _{i}-\partial _{i}\eta _{j})\xi ^{h}=0,
\label{eq-2}
\end{equation}%
where
\[
N_{ji}^{h}=\varphi _{j}^{t}\partial _{t}\varphi _{i}^{h}-\varphi
_{i}^{t}\partial _{t}\varphi _{j}^{h}-(\partial _{j}\varphi
_{i}^{t}-\partial _{i}\varphi _{j}^{t})\varphi _{t}^{h}
\]%
is the Nijenhuis tensor formed with $\varphi _{i}^{h}$ and $\partial
_{j}=\partial /\partial x^{j}$, then the almost contact structure is said to
be {\em normal} and the manifold is called a {\em normal almost contact
manifold}.

If, in an almost contact manifold, there is given a Riemannian metric $%
g_{ji} $ such that
\begin{equation}
g_{ts}\varphi _{j}^{t}\varphi _{i}^{s}=g_{ji}-\eta _{j}\eta _{i},\quad \eta
_{i}=g_{ih}\xi ^{h},\quad g_{kk}\xi ^{k}\xi ^{k}=1.  \label{eq-3}
\end{equation}
then the almost contact manifold is called an {\em almost contact metric
manifold }\cite{Blair,Sasaki-65}.

Comparing the first equations of (\ref{eq-1}) and (\ref{eq-3}), we see that
\[
\varphi _{ji}=\varphi _{j}^{t}g_{ti}
\]%
is skew-symmetric.

An almost contact metric manifold is said to be {\em cosymplectic manifold }%
\cite{Blair} if the $2$-form $\varphi _{ji}$ and the $1$-form $\eta _{j}$
are closed. It is well known that the cosymplectic manifold is characterized
by
\[
\nabla _{j}\varphi _{i}^{h}=0,\qquad \nabla _{i}\xi ^{h}=0,
\]%
where $\nabla _{j}$ denotes the operator of covariant differentiation with
respect to $g_{ji}$.

\subsection{Cosymplectic Bochner curvature tensor}

The cosymplectic Bochner curvature tensor \cite{Eum} is given by
\begin{eqnarray}
B_{kji}^{h} &=&K_{kji}^{h}+(\delta _{k}^{h}-\eta _{k}\eta
^{h})L_{ji}-(\delta _{j}^{h}-\eta _{j}\eta ^{h})L_{ki}  \nonumber \\
&&+L_{k}^{h}(g_{ji}-\eta _{j}\eta _{i})-L_{j}^{h}(g_{ki}-\eta _{k}\eta _{i})
\nonumber \\
&&+\varphi _{k}^{h}M_{ji}-\varphi _{j}^{h}M_{ki}+M_{k}^{h}\varphi
_{ji}-M_{j}^{h}\varphi _{ki}  \nonumber \\
&&-2(M_{kj}\varphi _{i}^{h}+\varphi _{kj}M_{i}^{h}),  \label{eq-Boc}
\end{eqnarray}%
where
\[
L_{ji}=-\frac{1}{2(m+2)}\left( K_{ji}+L(g_{ji}-\eta _{j}\eta _{i})\right)
,\quad M_{ji}=-L_{jt}\varphi _{i}^{t},
\]%
\[
L_{k}^{h}=L_{kt}g^{th},\quad M_{k}^{h}=M_{kt}g^{th},\quad L=g^{ji}L_{ji}=-%
\frac{K}{4(m+1)}.
\]%
By above equations, we can say that $L_{ji}$ is symmetric and $M_{ji}$ is
skew-symmetric.

It is easy to verify that
\[
L_{ji}\eta ^{i}=0,\quad M_{ji}\eta ^{i}=0.
\]%
The Bochner curvature tensor satisfies the follwing identities:%
\[
B_{kji}^{h}=-B_{jki}^{h},
\]%
\[
B_{kji}^{h}+B_{jik}^{h}+B_{ikj}^{h}=0,
\]%
\[
B_{hkji}=-B_{hkij},
\]%
\[
B_{hkji}=-B_{khji},
\]%
\[
B_{hkji}=B_{jihk},
\]%
\[
B_{kji}^{t}\eta _{t}=0,
\]%
\[
B_{kjti}\varphi _{i}^{t}=0,
\]%
where $B_{hkji}=B_{hkj}^{t}g_{ti}$.

\section{Cosymplectic conformal connections \label{sect-Co}}

Let $D$ be an affine connection in a cosymplectic manifold $M$ and $\Gamma
_{ji}^{h}$ be the components of the affine connection $D$. Assume that the
affine connection $D$ satisfies
\begin{equation}
D_{k}(e^{2p}g_{ji})=2e^{2p}p_{k}\eta _{j}\eta _{i},  \label{eq-wh-1}
\end{equation}%
where $D_{k}$ the operator of covariant differentiation with respect to $%
\Gamma _{ji}^{h}$ and $p$ is a certain scalar function, $p_{i}=\partial
_{i}p $. The torsion tensor of $D$ satisfies
\begin{equation}
\Gamma _{ji}^{h}-\Gamma _{ij}^{h}=-2\varphi _{ji}q^{h},  \label{eq-wh-2}
\end{equation}%
where $q^{h}$ is a certain vector field.

By solving (\ref{eq-wh-1}) and (\ref{eq-wh-2}), we get%
\begin{equation}
\Gamma _{ji}^{h}=\left\{
\begin{array}{c}
h \\
ij%
\end{array}%
\right\} +(\delta _{j}^{h}-\eta _{j}\eta ^{h})p_{i}+(\delta _{i}^{h}-\eta
_{i}\eta ^{h})p_{j}-(g_{ji}-\eta _{j}\eta _{i})p^{h}+\varphi
_{j}^{h}q_{i}+\varphi _{i}^{h}q_{j}-\varphi _{ji}q^{h},  \label{eq-wh-3}
\end{equation}%
where
\[
p^{h}=p_{t}g^{th},\qquad q^{h}=q_{t}g^{th}.
\]%
Using (\ref{eq-wh-3}), we can find that
\begin{equation}
D_{j}\varphi _{i}^{h}=(\delta _{j}^{h}-\eta _{j}\eta ^{h})(p_{t}\varphi
_{i}^{t}+q_{i})+(g_{ji}-\eta _{j}\eta _{i})(\varphi
_{t}^{h}p^{t}-q^{h})+\varphi _{j}^{h}(q_{t}\varphi _{i}^{t}-p_{i})+\varphi
_{ji}(p^{h}+\varphi _{t}^{h}q^{t}).  \label{eq-wh-4}
\end{equation}%
Now assume that affine connection satisfies
\begin{equation}
D_{j}\varphi _{i}^{h}=0.  \label{eq-wh-4.}
\end{equation}%
Using (\ref{eq-wh-4.}) in (\ref{eq-wh-4}) and contracting with respect to $h$
and $j$, we have
\begin{equation}
n(p_{t}\varphi _{i}^{t}+q_{i})+\eta _{i}\eta _{t}q^{t}=0,  \label{eq-w-1}
\end{equation}%
on transvecting with $\eta ^{i}$, we have
\begin{equation}
\eta _{t}q^{t}=0=q_{t}\eta ^{t}=0.  \label{eq-w-2}
\end{equation}%
Using (\ref{eq-w-2}) in (\ref{eq-w-1}), we get
\[
q_{i}=-p_{t}\varphi _{i}^{t}.
\]%
Now assume that the affine connection $D$ also satisfies
\begin{equation}
D_{j}\eta ^{h}=0.  \label{eq-wh-5}
\end{equation}%
So
\[
(\delta _{j}^{h}-\eta _{j}\eta ^{h})\eta _{t}p^{t}=0.
\]%
On contracting above equation with respect to $h$ and $j$, we obtain
\[
\eta _{t}p^{t}=0=\eta ^{t}p_{t}.
\]%
Therefore
\begin{equation}
p_{i}=q_{t}\varphi _{i}^{t},\quad q^{h}=\varphi _{t}^{h}p^{t},\quad
p^{h}=-\varphi _{t}^{h}q^{t}.  \label{eq-con-2}
\end{equation}%
We can easily calculate
\begin{equation}
p_{i}q^{i}=0,\quad p_{i}p^{i}=q_{i}q^{i}.  \label{eq-con-3}
\end{equation}%
By using (\ref{eq-wh-4}) and (\ref{eq-wh-5}), we obtain

\begin{prop}
Let $M$ be a cosymplectic manifold with structure tensors $(\varphi
_{i}^{h},\xi ^{h},\eta _{i},g_{ij})$ such that the affine connection $D$
satisfies
\[
D_{k}(e^{2p}g_{ji})=2e^{2p}p_{k}\eta _{j}\eta _{i},\qquad D_{j}\varphi
_{i}^{h}=0,\qquad D_{j}\eta ^{h}=0
\]%
and whose torsion tensor satisfies
\[
\Gamma _{ji}^{h}-\Gamma _{ij}^{h}=-2\varphi _{ji}q^{h},
\]%
where $p$ is a scalar function and $q^{h}$ is a vector field. Then the
components of the affine connection $D$ are given by
\begin{equation}
\Gamma _{ji}^{h}=\left\{
\begin{array}{c}
h \\
ji%
\end{array}%
\right\} +(\delta _{j}^{h}-\eta _{j}\eta ^{h})p_{i}+(\delta _{i}^{h}-\eta
_{i}\eta ^{h})p_{j}-(g_{ji}-\eta _{j}\eta _{i})p^{h}+\varphi
_{j}^{h}q_{i}+\varphi _{i}^{h}q_{j}-\varphi _{ji}q^{h},  \label{eq-cur}
\end{equation}%
where $p^{h}=p_{t}g^{th},u^{h}=u_{t}g^{th}$, $q_{i}=-p_{t}\varphi _{i}^{t}$,
$p_{i}=q_{t}\varphi _{i}^{t}.$
\end{prop}

We will call such an affine connection $D$ a {\em cosymplectic conformal
connection}.

\begin{prop}
A cosymplectic conformal connection satisfies
\[
D_{k}\left( e^{2p}(g_{ji}-\eta _{j}\eta _{i})\right) =0.
\]
\end{prop}

\section{Curvature tensor of a cosymplectic conformal connection\label%
{sect-curv}}

In this section, we will derive the curvature tensor of the cosymplectic
conformal connection.

Using (\ref{eq-cur}), the curvature tensor of the cosymplectic conformal
connection is
\begin{eqnarray}
R_{kji}^{h} &=&K_{kji}^{h}-(\delta _{k}^{h}-\eta _{k}\eta
^{h})p_{ji}+(\delta _{j}^{h}-\eta _{j}\eta ^{h})p_{ki}  \nonumber \\
&&-p_{k}^{h}(g_{ji}-\eta _{j}\eta _{i})+p_{j}^{h}(g_{ki}-\eta _{k}\eta _{i})
\nonumber \\
&&-\varphi _{k}^{h}q_{ji}+\varphi _{j}^{h}q_{ki}-q_{k}^{h}\varphi
_{ji}+q_{j}^{h}\varphi _{ki}  \nonumber \\
&&+(\nabla _{k}q_{j}-\nabla _{j}q_{k})\varphi _{i}^{h}+2\varphi
_{kj}(q_{i}p^{h}-p_{i}q^{h}),  \label{eq-R}
\end{eqnarray}%
where
\begin{equation}
p_{ji}=\nabla _{j}p_{i}-p_{j}p_{i}+q_{j}q_{i}+\frac{1}{2}\lambda
(g_{ji}-\eta _{j}\eta _{i}),  \label{eq-p1}
\end{equation}
\begin{equation}
q_{ji}=\nabla _{j}q_{i}-p_{j}q_{i}-p_{i}q_{j}+\frac{1}{2}\lambda \varphi
_{ji},  \label{eq-p2}
\end{equation}%
\[
p_{k}^{h}=p_{kt}g^{th},\quad q_{k}^{h}=q_{kt}g^{th},\quad \beta
_{k}^{h}=\beta _{kt}g^{th}.
\]%
$\lambda $ being defined by
\[
\lambda =p_{i}p^{i}=q_{i}q^{i}.
\]%
Let
\[
\alpha _{ji}=-(\nabla _{j}q_{i}-\nabla _{i}q_{j}),
\]%
\[
\beta _{ji}=2(p_{j}q_{i}-p_{i}q_{j}).
\]%
Since $p_{i}=\partial _{i}p$, so $p_{ji}=p_{ij}$. We can easily check that
\begin{equation}
\eta ^{j}p_{ji}=0,\quad \eta ^{j}q_{ji}=0,\quad \alpha _{ji}\eta ^{i}=0,
\end{equation}%
\begin{equation}
\beta _{ji}\eta ^{i}=0,\quad \alpha _{ji}=-\alpha _{ij},\quad \beta
_{ji}=-\beta _{ij},
\end{equation}%
\begin{equation}
q_{ji}=-p_{jt}\varphi _{i}^{t},\quad p_{ji}=q_{jt}\varphi _{i}^{t},
\label{eq-cc}
\end{equation}%
\begin{equation}
\alpha _{ji}=-(q_{ji}-q_{ij}-\lambda \varphi _{ji}),  \label{eq-q3}
\end{equation}%
\begin{equation}
\alpha =\varphi ^{ij}\alpha _{ij}=-2\nabla _{t}p^{t},  \label{eq-q1}
\end{equation}%
\begin{equation}
\beta =\varphi ^{ij}\beta _{ij}=4\lambda .  \label{eq-q2}
\end{equation}%
By (\ref{eq-p1}), we have
\begin{equation}
p_{k}^{k}=\nabla _{k}p^{k}+mp_{k}p^{k}.  \label{eq-p3}
\end{equation}%
We now consider here that the curvature tensor of cosymplectic curvature
connection vanishes, that is,
\begin{equation}
R_{kji}^{h}=0.  \label{eq-11}
\end{equation}%
Consequently (\ref{eq-R}) becomes
\begin{eqnarray}
K_{kji}^{h} &=&(\delta _{k}^{h}-\eta _{k}\eta ^{h})p_{ji}-(\delta
_{j}^{h}-\eta _{j}\eta ^{h})p_{ki}  \nonumber \\
&&+p_{k}^{h}(g_{ji}-\eta _{j}\eta _{i})-p_{j}^{h}(g_{ki}-\eta _{k}\eta _{i})
\nonumber \\
&&+\varphi _{k}^{h}q_{ji}-\varphi _{j}^{h}q_{ki}+q_{k}^{h}\varphi
_{ji}-q_{j}^{h}\varphi _{ki}  \nonumber \\
&&+\alpha _{kj}\varphi _{i}^{h}-\varphi _{kj}\beta _{i}^{h},  \nonumber
\end{eqnarray}%
in covariant form,
\begin{eqnarray}
K_{kjih} &=&(g_{kh}-\eta _{k}\eta _{h})p_{ji}-(g_{jh}-\eta _{j}\eta
_{h})p_{ki}  \nonumber \\
&&+p_{kh}(g_{ji}-\eta _{j}\eta _{i})-p_{jh}(g_{ki}-\eta _{k}\eta _{i})
\nonumber \\
&&+\varphi _{kh}q_{ji}-\varphi _{jh}q_{ki}+q_{kh}\varphi _{ji}-q_{jh}\varphi
_{ki}  \nonumber \\
&&+\alpha _{kj}\varphi _{ih}-\varphi _{kj}\beta _{ih},  \label{eq-ct}
\end{eqnarray}%
Using (\ref{eq-ct}) in the identity $K_{kjih}=K_{ihkj}$, we have
\begin{eqnarray}
0 &=&\varphi _{ji}\left( q_{kh}+q_{hk}\right) -\varphi _{jh}\left(
q_{ki}+q_{ik}\right) +\varphi _{kh}\left( q_{ji}+q_{ij}\right)  \nonumber \\
&&-\varphi _{ki}\left( q_{jh}+q_{hj}\right) +\alpha _{kj}\varphi
_{ih}-\alpha _{ih}\varphi _{kj}-\beta _{ih}\varphi _{kj}+\beta _{kj}\varphi
_{ih},  \label{eq-qq}
\end{eqnarray}%
transvecting (\ref{eq-qq}) with $\varphi ^{kh}$, we get
\begin{equation}
q_{ji}+q_{ij}=0.  \label{eq-qij}
\end{equation}%
Using (\ref{eq-qij}) in (\ref{eq-qq}), we obtain
\[
\left( \alpha _{kj}+\beta _{kj}\right) \varphi _{ih}=\left( \alpha
_{ih}+\beta _{ih}\right) \varphi _{kj},
\]%
transvecting the above equation with $\varphi ^{ih}$, we get
\[
\alpha _{kj}+\beta _{kj}=\frac{1}{2m}\left( \alpha +\beta \right) \varphi
_{kj},
\]%
using (\ref{eq-q1}), (\ref{eq-q2}), we have
\[
\alpha _{kj}+\beta _{kj}=\frac{1}{m}\left( -\nabla
_{t}p^{t}+2p_{t}p^{t}\right) \varphi _{kj}.
\]%
Using (\ref{eq-q3}) and (\ref{eq-qij}) in above equation, we get
\[
\beta _{kj}=\frac{1}{m}\left( -\nabla _{t}p^{t}+(2-m)p_{t}p^{t}\right)
\varphi _{kj}+2q_{kj}.
\]%
By use of (\ref{eq-p3}) in above equation, we have
\begin{equation}
\beta _{kj}=\frac{1}{m}\left( -p_{t}^{t}+2p_{t}p^{t}\right) \varphi
_{kj}+2q_{kj}.  \label{eq-B}
\end{equation}%
By use of (\ref{eq-ct}) in the identity $K_{kjih}+K_{jikh}+K_{ikjh}-0$ and
using~(\ref{eq-q3}), (\ref{eq-qij}), (\ref{eq-B}), we have%
\[
\frac{1}{m}\left( (m-2)p_{t}p^{t}+p_{t}^{t}\right) \left( \varphi
_{kj}\varphi _{ih}+\varphi _{ji}\varphi _{kh}+\varphi _{ik}\varphi
_{jh}\right) =0,
\]%
from which, we obtain
\begin{equation}
(m-2)p_{t}p^{t}+p_{t}^{t}=0.  \label{eq-BB}
\end{equation}%
Substituting (\ref{eq-B}) in (\ref{eq-BB}), we have
\begin{equation}
\beta _{kj}=p_{t}p^{t}\varphi _{kj}+2q_{kj}.  \label{eq-BBB}
\end{equation}%
By (\ref{eq-cc}), we have
\[
q_{ji}\varphi _{s}^{i}\varphi _{k}^{j}=-q_{sk}.
\]%
Using (\ref{eq-BBB}), we get
\begin{equation}
\beta _{jk}\varphi _{s}^{j}=-2p_{ks}-p_{t}p^{t}(g_{ks}-\eta _{k}\eta _{s}).
\label{eq-A}
\end{equation}%
By use of (\ref{eq-q3}) and (\ref{eq-qij}), we have
\begin{equation}
\alpha _{kj}\varphi _{s}^{j}=2p_{ks}-p_{t}p^{t}(g_{ks}-\eta _{k}\eta _{s}).
\label{eq-AA}
\end{equation}%
Contracting (\ref{eq-ct}) with respect to $h$ and $k$ and using (\ref{eq-A}%
), (\ref{eq-AA}), we obtain
\begin{equation}
K_{ji}=(2m+4)p_{ji}+p_{t}{}^{t}(g_{ji}-\eta _{j}\eta _{i}),  \label{eq-ct1}
\end{equation}%
Transvecting (\ref{eq-ct1}) with $g^{ji}$, we find
\begin{equation}
K=(4m+4)p_{t}{}^{t}.  \label{eq-ct4}
\end{equation}%
Consequently, we have
\begin{equation}
p_{t}{}^{t}=\frac{K}{4(m+1)}=-L  \label{eq-a}
\end{equation}%
Using (\ref{eq-a}) in (\ref{eq-ct1}), we have
\[
K_{ji}=(2m+4)p_{ji}-L(g_{ji}-\eta _{j}\eta _{i}),
\]%
so
\begin{equation}
p_{ji}=\frac{1}{(2m+4)}\left( K_{ji}+L(g_{ji}-\eta _{j}\eta _{i})\right)
=-L_{ji}.  \label{eq-L}
\end{equation}%
Therefore
\begin{equation}
M_{ji}=-L_{jt}\varphi _{i}^{t}=p_{jt}\varphi _{i}^{t}=-q_{ji}.  \label{eq-M}
\end{equation}%
From (\ref{eq-BB}), we have
\[
p_{t}p^{t}=\frac{L}{(m-2)}.
\]%
By use of (\ref{eq-BBB}), (\ref{eq-q3}) and (\ref{eq-qij}), we obtain
\begin{equation}
\alpha _{ji}=2M_{ji}+\frac{L}{(m-2)}\varphi _{ji},  \label{eq-LL}
\end{equation}%
\begin{equation}
\beta _{ji}=-2M_{ji}+\frac{L}{(m-2)}\varphi _{ji}.  \label{eq-LLL}
\end{equation}%
Substituting (\ref{eq-ct}), (\ref{eq-L}), (\ref{eq-M}), (\ref{eq-LL}), (\ref%
{eq-LLL}) in (\ref{eq-Boc}), we obtain
\[
B_{kjih}=0.
\]%
Thus, we have the following result:

\begin{th}
Let $M$ be a real $(2m+1)$-dimensional $(m>2)$ cosymplectic manifold. If $M$
admits a cosymplectic conformal connection which is of zero curvature, then
the Bochner curvature tensor of $M$ vanishes.
\end{th}

Department of Mathematics \& Statistics

School of Mathematical \& Physical Sciences

Dr. Harisingh Gour University

Sagar - 470 003, MP, INDIA

Email: pgupta@dhsgsu.edu.in

\end{document}